
\def\a{\alpha}
\def\b{\beta}
\def\d{\delta}
\def\D{\Delta}
\def\g{\gamma}
\def\G{\Gamma}
\def\lb{\lambda}
\def\Om{\Omega}

\def\vp{\varphi}
\def\z{\zeta}

\font\tenbb=msym10
\font\sevenbb=msym7
\font\fivebb=msym5
\newfam\bbfam
\textfont\bbfam=\tenbb \scriptfont\bbfam=\sevenbb
\scriptscriptfont\bbfam=\fivebb
\def\bb{\fam\bbfam}

\def\Ab{{\bb A}}
\def\Cb{{\bb C}}
\def\Pb{{\bb P}}
\def\Rb{{\bb R}}

\def\Oc{{\cal O}}

\def\part{\partial}
\def\ov{\overline}
\def\ify{\infty}
\def\lgl{\langle}
\def\rgl{\rangle}
\def\ra{\rightarrow}
\def\sbs{\subset}
\def\ot{\otimes}
\def\wt{\widetilde}
\def\nb{\nabla}

\def\tr{\mathop{\rm tr}\nolimits}
\def\dt{\mathop{\rm dt}\nolimits}

\def\ch{\mathop{\rm ch}\nolimits}
\def\Td{\mathop{\rm Td}\nolimits}
\def\Ker{\mathop{\rm Ker}\nolimits}

\catcode`\@=11
\def\displaylinesno #1{\displ@y\halign{
\hbox to\displaywidth{$\@lign\hfil\displaystyle##\hfil$}&
\llap{$##$}\crcr#1\crcr}}

\def\ldisplaylinesno #1{\displ@y\halign{
\hbox to\displaywidth{$\@lign\hfil\displaystyle##\hfil$}&
\kern-\displaywidth\rlap{$##$}
\tabskip\displaywidth\crcr#1\crcr}}
\catcode`\@=12

\def\build#1_#2^#3{\mathrel{
\mathop{\kern 0pt#1}\limits_{#2}^{#3}}}

\magnification=1200
\overfullrule=0pt

\hsize=125mm
\vsize=180mm
\hoffset=3mm
\voffset=12mm

\baselineskip=14pt

\def\ify{\infty}
\def\lgl{\langle}
\def\nb{\nabla}
\def\op{\oplus}
\def\ot{\otimes}
\def\part{\partial}
\def\rgl{\rangle}
\def\sbs{\subset}
\def\sm{\simeq}

\def\a{\alpha}
\def\b{\beta}
\def\lb{\lambda}
\def\t{\theta}
\def\vp{\varphi}
\def\z{\zeta}

\def\D{\Delta}
\def\L{\Lambda}

\def\Ker{\mathop{\rm Ker}\nolimits}
\def\rank{\mathop{\rm rank}\nolimits}
\def\Re{\mathop{\rm Re}\nolimits}
\def\Tr{\mathop{\rm Tr}\nolimits}

\def\ra{\rightarrow}
\def\longra{\longrightarrow}
 
\def\build#1_#2^#3{\mathrel{
\mathop{\kern 0pt#1}\limits_{#2}^{#3}}}
 
\def\Ec{{\cal E}}
\def\Oc{{\cal O}}
\def\Xc{{\cal X}}
 
\font\tenbb=msym10
\font\sevenbb=msym7
\font\fivebb=msym5
\newfam\bbfam
\textfont\bbfam=\tenbb \scriptfont\bbfam=\sevenbb
\scriptscriptfont\bbfam=\fivebb
\def\bb{\fam\bbfam}
 
\def\Cb{{\bb C}}
\def\Pb{{\bb P}}
\def\Rb{{\bb R}}
\def\Zb{{\bb Z}}

\centerline{\bf Upper bounds for regularized determinants}

\bigskip

\centerline{by}

\medskip

\centerline{H. Gillet \footnote*{ Supported
by N.S.F grant  DMS-9501500} and C. Soul\'e}

\vglue 1cm

Let $E$ be a holomorphic vector bundle on a compact
K\"ahler manifold $X$. If we fix a metric $h$ on $E$, we
get a Laplace operator $\D$ acting upon smooth sections
of $E$ over $X$. Using the zeta function of $\D$, one
defines its regularized determinant $\det' (\D)$. In [5]
\S4.1.6, inspired by our arithmetic Riemann-Roch theorem,
we were led to conjecture that, when $h$ varies, this
determinant $\det' (\D)$ remains bounded from above.

In this paper we prove this in two special cases. The
first case is when $X$ is a Riemann surface, $E$ is a
line bundle and $\dim H^0 (X,E) + \dim H^1 (X,E) \leq
2$, and the second case is when $X=\Pb^1$, $E$ is a line
bundle, and all metrics under consideration
are invariant under rotation around a fixed axis.
To get the desired upper bound in the first case we use
an inequality of Moser and Trudinger, and its
extension to arbitrary compact manifolds due to Fontana
[3]. We prove the second case by direct estimates.

Though our results deal with very few cases, we 
find striking that inequalities as sharp as the theorem
of Moser and Trudinger can be used to prove
our conjecture. We hope
the reader will get interested in the general question,
and try to either prove or disprove our statement.

\medskip

In the first paragraph we phrase the conjecture
in its most general form, 
and give a few facts about it.
Next, in the case of a line bundle over a Riemann surface,
 we compute the anomaly $A(\vp)$
for the regularized
determinant of the Laplace operator
when a fixed metric $h_0$ on the line bundle 
is replaced by $h_0 e^{\vp}$.
To check our theorem, we then need to 
bound from above the functional
$A(\vp)$ when $\vp$
is any smooth function on the Riemann surface
(resp. any function of the distance to the origin on the
projective line).
This is done in the next two 
paragraphs.
At the end, we discuss the (much easier)
case of the trivial line bundle on the circle, 
where the determinant is bounded from below.
\medskip

We thank  W. Beckner, P. Chang and
J. Lott for interesting discussions.

\vglue 1cm

\noindent {\bf 1. Statement of the results.}

\medskip
\noindent {\bf 1.1.} Let $X$ be a smooth, projective,
equidimensional complex variety of dimension $d$
 and $h_X$ an hermitian metric on
its tangent space. The associated (normalized) K\"ahler
form $\mu$ is defined by the formula
$$
\mu = {i\over 2\pi}\  \sum_{\a,\b}\ h_X
 \left( {\part \over {\part z}_{\a}} \,
, \, {\part \over {\part z}_{\b}}\right) \, dz_{\a} \, d\ov z_{\b} \, ,
$$
where $(z_{\a})$ is any local holomorphic chart on $X$.

\medskip

Consider an holomorphic vector bundle $E$ on $X$, equipped
with a $C^{\ify}$ hermitian metric $h$. Let $A^{0q}
(X,E)$, $q=0,..,d$, be the space of smooth forms of type
$(0,q)$ with values in $E$. The $L^2$-metric on $A^{0q}
(X,E)$ is defined by the formula
$$
\lgl s,t \rgl_{L^2} = \int_X \lgl s(x),t(x)\rgl \, \mu^d/d! \,
,
$$
where $s,t \in A^{0q} (X,E)$ and $\lgl s(x),t(x)\rgl$ is
the pointwise scalar product defined by $h$ and $h_X$
([5], \S4.1.1). The Cauchy-Riemann operator
$$
\ov{\part} :\ A^{0q} (X,E) \ra A^{0,q+1} (X,E)
$$
has an adjoint $\ov{\part}^*$:
$$
\lgl \ov{\part} \, s , t\rgl_{L^2} = \lgl s,\ov{\part}^* \,
t\rgl_{L^2} \, .
$$
We consider the Laplace operator $\D_q = \ov{\part}^* \,
\ov{\part}$ on $A^{0q} (X,E)$ and its zeta function
$$
\z_{\D_q} (s) = {1\over \G (s)} \int_0^{\ify} \tr
 \, (e^{-\D_q
t}) \, t^{s-1} \, \dt \ , \ \Re \, (s) > 1 \, .
$$
It is known that $\z_{\D_q} (s)$ can be analytically
continued to the whole complex plane and is regular at
the origin. The regularized determinant of $\D_q$ is
defined to be
$$
{\det}' (\D_q) = \exp \, (-\z'_{\D_q} (0)) \, ,
$$
where $\z'_{\D_q} (0)$ is the value of ${d \over ds} \,
\z_{\D} (s)$ at $s=0$. Our goal is to find upper bounds
for ${\det}' (\D_q)$ when $h$ varies.

More precisely, for any $q \geq 0$, consider the spaces
$$
B^q = \bar{\part} (A^{0,q-1} (X,E)) \sbs A^{0q} (X,E) \, , \quad
q \geq 1 \, ,
$$
$B^0 =0$, and the zeta function
$$
\z_{B^q} (s) = \Tr (\D_q^{-s} \mid B^q) \, , \quad \Re (s) > d
\, .
$$
By the Hodge decomposition theorem we have
$$
A^{0q} (X,E) = B^q \op \bar{\part}^* (A^{0,q+1} (X,E)) \op \Ker
(\D_q) \, ,
$$
and the Cauchy-Riemann operator induces an isomorphism
$$
\bar{\part} : \bar{\part}^* (A^{0q} (X,E)) \build
\longra_{}^{\sim} B^q
$$
such that $\bar{\part} \, \D_{q-1} = \D_q \, \bar{\part}$. It
follows that
$$
\z_{\D_q} (s) = \z_{B^q} (s) + \z_{B^{q+1}} (s)
$$
hence
$$
\z_{B^{q+1}} (s) = \z_{\D_q} (s) - \z_{\D_{q-1}} (s) +
\z_{\D_{q-2}} (s) + \cdots + (-1)^q \, \z_{\D_0} (s) \, .
$$
This implies that $\z_{B^q} (s)$ converges when $\Re (s) > d$,
has a meromorphic continuation to the whole complex plane, and
is regular at the origin. Define
$$
D_q (E,h) = \exp (-\z'_{B^q} (0)) \, .
$$
In [5], \S 4.1.6, we proposed the following

\bigskip

\noindent {\bf Conjecture.} {\it There exists a constant $C_q (E)$
such that,
for any choice of a metric $ h$ on $ E$}  , 
$$D_q (E,h) \leq C_q (E) \, .$$

\bigskip

\noindent {\bf 1.2. Remarks on the conjecture}

\smallskip

\noindent {\bf 1.2.1.} In the conjecture above both $D_q (E,h)$
and $C_q (E)$ depend in general on the metric $h_X$ on $X$.

Notice that, for any real constant $t>0$,
$$
D_q (E,th) = D_q (E,h) \, . \leqno (1)
$$
Indeed, when $h$ is replaced by $th$, the $L^2$-metric on
$A^{0q} (X,E)$ gets multiplied by 
the same factor $t$ for all $q \geq 0$,
therefore $\D_q$ remains unchanged.

Furthermore, if $E^\vee = E^* \ot \L^d (TX^*)$ is  Poincar\'e
dual of $E$ and $h^\vee$ the metric on $E^\vee$ induced by $h$
and $h_X$, the Poincar\'e-Serre duality implies that
$$
D_q (E,h) = D_{d+1-q} (E^\vee ,h^\vee) \, . \leqno (2)
$$

Therefore the conjecture is stable under scaling and duality.

\medskip

\noindent {\bf 1.2.2.} Our inspiration to make this conjecture
was number theoretic. Assume that $X$ is the set $\Xc (\Cb)$ of
complex points of a regular projective flat scheme $\Xc$ over
$\Zb$ and that $E$ is the holomorphic vector bundle defined by
an algebraic vector bundle $\Ec$ on $\Xc$. In loc.cit. we defined
arithmetic Betti numbers
$$
b_q (\Ec ,h) \in \Rb \ , \quad 0 \leq q \leq d+1 \, ,
$$
as follows .
When $M$ is a finitely generated abelian group, equipped with a
norm $\Vert \cdot \Vert$ on its real span $M \, \build
\ot_{\Zb}^{} \, \Rb$, we let
$$
h^0 (M, \Vert \cdot \Vert) = \log \, \# \, \{ m \in M / \Vert m
\Vert \leq 1 \}
$$
and
$$
h^1 (M, \Vert \cdot \Vert ) = h^0 (M^* , \Vert \cdot \Vert^*) \,
.
$$
We consider the  coherent cohomology groups
$H^q (\Xc ,\Ec)$,
$q \geq 0$, equipped with their $L^2$ metric. Then we let
$$
b_q (\Ec ,h) = h^0 (H^q (\Xc ,\Ec)) + h^1 (H^{q-1} (\Xc ,\Ec)) +
{1\over 2} \, \z'_{B^q} (0) \, . \leqno (3)
$$
In [5], loc.cit., we gave properties of these numbers (duality,
Euler characteristic formula) which partially justified
calling
them Betti numbers. However, a basic property should be that
each $b_q (\Ec ,h)$ is nonnegative, or at least bounded below.
This led us to the conjecture in \S1.1.

\medskip

\noindent {\bf 1.2.3.} It would be of interest to find some
 interpretation of our conjecture in mathematical physics.
In [6], the Moser-Trudinger inequality
(see  \S 2.3 below) is interpreted as the existence of
a lower bound for a free energy
functional, and it is derived
in op.cit. Prop. 4 from the
 Gibbs variational principle.

\medskip

\noindent {\bf 1.2.4.} A stronger version of the conjecture
consists in requiring that $C_q (E)$ depends only on the
$C^{\ify}$ bundle underlying $E$, and not on its holomorphic
structure. Results like [10] Proposition 3 (due to Miyaoka and
based on a result of Selberg, [10] Theorem 4) , which says
that when $d=1$ and when $E$ is a flat unitary bundle,
the following holds
$$
\log D_1 (E,h) \leq  constant \cdot  \rank (E) \, ,
$$
points in this direction.

\bigskip

\noindent {\bf 1.3. Results}

\smallskip

From now on we assume that $X$ is a curve $(d=1)$ and that $E$
is a line bundle $L$. We then take $q=1$ and we write $D(L,h)$
instead of $D_1 (E,h) = \det' (\D_0)$.

\bigskip

\noindent {\bf Theorem 1.}

\smallskip

\item{i)}{\it If
$$
\dim H^0 (X,L) + \dim H^1 (X,L) \leq 2 \, ,
$$
there exists a constant $C(L)$ such that, for any metric $h$ on
$L$,
$$
D(L,h) \leq C(L) \, ;
$$

\item{ii)} Assume that $X = \Pb^1 (\Cb)$ is the complex
projective line, that $h_X$ is invariant under rotation, and
that $L=\Oc (n)$, $n \in \Zb$. Then there is a constant $C(n)$
such that, for any metric $h$ on $L$ invariant under rotation,}
$$
D(L,h) \leq C(n) \, .
$$

\bigskip

To clarify statement ii) above, let us write $z=r \, e^{i\t}$
the standard coordinate on $\Cb \sbs \Pb^1 (\Cb)$. Given any
$\a$, we let $r_{\a} (z) = e^{i\a} \, z$ be the rotation of
angle $\a$. A metric $h_X$ on $X$ (resp. $h$ on $L$) is said to
be invariant under rotation when $r_{\a}^* (h_X) = h_X$ (resp.
$r_{\a}^* (h) = h$) for all values of $\a$.

\bigskip

\noindent {\bf 2. An anomaly formula.}

\medskip

\noindent {\bf 2.1.} We fix $X$, $h_X$ and $L$ as in
\S 1.3. Let $h_0$ and $h=h_0 \exp \, (\vp)$ be two
hermitian metrics on $L$, with $\vp$ a smooth real
valued function on $X$. 
We shall give a formula comparing the determinants
$D(L,h)$ and
$D(L,h_0)$.

\medskip

Let $b_0 = \dim H^0 (X,L)$ and $b_1 =\dim H^1 (X,L)$. We
endow $H^0 (X,L) = \ker \, (\ov{\part}) \sbs A^{00} (X,L)$
and $H^1 (X,L) = \ker \, (\ov{\part}^*) \sbs A^{01} (X,L)$
with the $L^2$-metric coming from $h_0$. Let $(\a_i)$,
$i=1,\ldots ,b_0$, be an orthonormal basis of $H^0
(X,L)$, and $(\b_i)$, $i=1,\ldots ,b_1$, an orthonormal
basis of $H^1 (X,L)$. If $1\leq i,j \leq b_0$ we let
$\lgl \a_i ,\a_j \rgl$ be the pointwise scalar product of
$\a_i$ with $\a_j$. We define similarly $\lgl \b_i ,\b_j
\rgl$, $1\leq i,j \leq b_1$.

\medskip

If $d=\part + \ov{\part}$ we let $d^c = {\part - \ov{\part}
\over 4\pi i}$, so that $dd^c = {\ov{\part} \, \part
\over 2\pi i}$.

\medskip

Denote by $c_1 (T_X ,h_X)$ the first Chern form of the
tangent bundle to $X$, and by $c_1 (L,h_0)$ the first
Chern form of $(L,h_0)$.
Clearly, to prove Theorem 1 it will be enough to show that,
under the given hypotheses, the quantity
$$
A(\vp) = \log D (L,h) - \log D (L,h_0) \leqno (4)
$$
remains bounded from above when $\vp$ varies.

\bigskip

\noindent {\bf Proposition 1.} {\it The following
formula holds:}
$$
\eqalign{
& \ A(\vp) \cr
= & \ {1\over 2} \int_X \vp \, dd^c (\vp) - \int_X \vp
(c_1 (T_X ,h_X)+c_1 (L,h_0)) \cr
+ & \ \log \det \left( \int_X e^{\vp} \, \lgl \a_i ,\a_j
\rgl \, \mu \right)_{1\leq i,j \leq b_0} \cr
+ & \ \log \det \left( \int_X e^{-\vp} \, \lgl \b_i ,\b_j
\rgl \, \mu \right)_{1\leq i,j\leq b_1} \, . \cr
}
$$

\bigskip

\noindent {\bf 2.2.} To prove Proposition 1 we consider
the Quillen metric [9] on the complex line
$$
\lb (L) = \L^{b_0} \, H^0 (X,L) \ot (\L^{b_1} \, H^0
(X,L))^* \, .
$$
It is defined as the quotient of the $L^2$-metric by the
determinant of the Laplace operator:
$$
\Vert \cdot \Vert_{Q,h}^2 = \Vert \cdot \Vert_{L^2 ,h}^2
\ D(L,h)^{-1} \, .
$$

Therefore we get
$$
 \ A(\vp) 
=  \ \log \, {\Vert \cdot \Vert_{L^2 ,h}^2 \over \Vert
\cdot \Vert_{L^2 ,h_0}^2} - \log \, {\Vert \cdot
\Vert_{Q,h}^2 \over \Vert \cdot \Vert_{Q,h_0}^2} \, , \leqno (5) 
$$
and we are led to compute the variation of both the
$L^2$-metric and the Quillen metric on $\lb (L)$.

\medskip

Concerning the Quillen metric, we can use the anomaly
formula in [1], Theorem 0.3. Let $\wt{\ch} \, (h,h_0)$ be the
Bott-Chern secondary characteristic class of $L$, which
satisfies 
$$
dd^c \, \wt{\ch} \, (h,h_0) = \ch \, (L,h) - \ch \,
(L,h_0) \, , 
$$
where $\ch \, (L,h) = \exp \, (c_1 (L,h))$ is the Chern
character form of $(L,h)$. If $\Td \, (T_X ,$ $h_X)$ is the
form representing  the Todd class of $X$ defined using
$h_X$, the following holds (loc.cit.):
$$
-\log \, {\Vert \cdot \Vert_{Q,h}^2 \over \Vert \cdot
\Vert_{Q,h_0}^2} = \int_X \wt{\ch} \, (h,h_0) \, \Td
\, (T_X ,h_X) \, .
$$
Since $X$ has dimension one, we have
$$
\eqalign{
& \ \int_X \wt{\ch} \, (h,h_0) \, \Td \, (T_X ,h_X) \cr
= & \ \int_X \wt{\ch}_2 (h,h_0) + \int_X \wt{\ch}_1
(h,h_0) \, c_1 (T_X ,h_X) \, , \cr
}
$$
where $\wt{\ch}_p$ is the component of degree $(p-1,p-1)$
of $\wt{\ch}$.

\medskip

Let $(\Oc_X ,e^{\vp})$ be the trivial line bundle
equipped with the metric such that $\Vert 1 \Vert^2 =
\exp \, (\vp)$. Since
$$
(L,h) = (L,h_0) \ot (\Oc_X ,e^{\vp})
$$
we deduce from [4], Proposition 1.3.3 (and formula
(1.3.5.2)) that
$$
\wt{\ch} \, (h,h_0) = \ch \, (L,h_0) \, \wt{\ch}
\, (e^{\vp}, 1) \, .
$$
To compute the Bott-Chern class $\wt{\ch} \, (e^{\vp} ,1)$
comparing the metric such that $\Vert 1 \Vert^2 = \exp \,
(\vp)$ with the trivial metric on $\Oc_X$, first notice
that $$
\wt{\ch}_1 (e^{\vp} ,1) = \wt{c}_1 (e^{\vp} ,1) =-\vp
$$
by [4], (1.2.5.1). Furthermore, from [4] (1.3.1.2) we get
$$
\displaylines{
\wt{\ch}_2 (e^{\vp} ,1) = {1\over 2} \, \wt{c}_1^2
(e^{\vp} ,1) \cr
= {1\over 2} \, c_1 (e^{\vp}) \, \wt{c}_1 (e^{\vp} ,1) =
{1\over 2} \, dd^c (\vp) \, \vp \, . \cr
}
$$
So we conclude that
$$
-\log \, {\Vert \cdot \Vert_{Q,h}^2 \over \Vert \cdot
\Vert_{Q,h_0}^2} =  \ -\int_X c_1 (L,h_0) \, \vp +
{1\over 2} \int_X \vp \, dd^c (\vp)
 \ - \int_X \vp \, c_1 (T_X ,h_X) \, . \leqno (6)
$$

\bigskip

\noindent {\bf 2.3.} Now we have to compute the variation
of the $L^2$-norm on $\lb (L)$. Since $(\a_i)$ is an
orthonormal basis of $\Ker \, (\ov{\part})$ for the
$L^2$-metric defined by $h_0$, the change of metric on
$\L^{b_0} \, H^0 (X,L)$ is the determinant
$$
\leqalignno{
q_0 = & \ \det \, (\lgl \a_i ,\a_j \rgl_{L^2 ,h}) &(7) \cr
= & \ \det \left( \int_X e^{\vp} \lgl \a_i ,\a_j \rgl \,
\mu \right) \, . \cr
}
$$

\medskip

For $H^1$ the situation is more complicated since the
kernel of $\ov{\part}^*$ changes when $h_0$ is replaced
by $h$. If $\psi_i \in A^{01} (X,L)$ is the harmonic form
for the metric \break coming from $h$ which is cohomologous
to $\b_i$, $1\leq i \leq b_1$, the change of metric on
$\L^{b_1} \, H^1 (X,L)$ is
$$
q_1 = \det (\lgl \psi_i ,\psi_j \rgl_{L^2 ,h}) \, .
\leqno (8)
$$
Note that $\ov{\part}_h^* (e^{-\vp} \, \b_i)=0$ since,
for any section $s\in A^{00} (L)$,
$$
\int_X \lgl \ov{\part} \, s , e^{-\vp} \, \b_i \rgl_h \,
\mu = \int_X \lgl \ov{\part} \, s , \b_i \rgl \, \mu =0
\, .
$$
Therefore there are complex numbers $a_{ij}$, $1\leq i,j
\leq b_1$, such that
$$
e^{-\vp} \, \b_i = \sum_{j=1}^{b_1} a_{ij} \, \psi_j \, .
\leqno (9)
$$
It follows that
$$
\det \, (\lgl \psi_i ,\psi_j \rgl_{L^2 ,h}) = \det \,
(\lgl e^{-\vp} \, \b_i , e^{-\vp} \, \b_j \rgl_{L^2 ,h}) \,
\det \, (a_{ij})^{-2} \, . \leqno (10)
$$
Note that
$$
\lgl e^{-\vp} \, \b_i , e^{-\vp} \, \b_j \rgl_{L^2 ,h} =
\int_X e^{-\vp} \, \lgl \b_i ,\b_j \rgl \, \mu \, .
\leqno (11)
$$
To compute $(a_{ij})$ let us introduce $\g_i = * \,
\b_i$, $1\leq i \leq b_1$, where $*$ is the star operator
defined by $h_X$ and $h_0$. Then $\g_i$ is a smooth
section of the Serre dual $\Om_X^1 \ot L^*$ of $L$ such
that $\ov{\part} \, (\g_i) =0$ (since $\ov{\part}_{h_0}^*
(\b_i) =0$) and
$$
\lgl \b_i , \b_j \rgl_{h_0} \, \mu = \b_i \, \g_j \quad ,
\quad 1\leq i,j \leq b_1 \, .
$$
Since $(\b_i)$ is orthonormal and $\psi_i$ is
cohomologous to $\b_i$ we get
$$
\int_X \psi_i \, \g_j = \int_X \b_i \, \g_j = \d_{ij} \, .
$$
Using (9) we deduce that
$$
\int_X e^{-\vp} \, \b_i \, \g_k = \sum_{j=1}^{b_1} a_{ij}
\int_X \psi_j \, \g_k = a_{ik} \, .
$$
In other words
$$
a_{ik} = \int_X e^{-\vp} \, \lgl \b_i ,\b_k\rgl_{h_0} \,
\mu \, . \leqno (12)
$$
From (8), (10), (11), (12) we conclude that
$$
q_1 = \det \left( \int_X e^{-\vp} \, \lgl \b_i ,\b_j
\rgl_{h_0} \, \mu \right)^{-1} \, . \leqno (13)
$$
From (7) and (13) we get
$$
\ldisplaylinesno{
\log \, {\Vert \cdot \Vert_{L^2 ,h}^2 \over \Vert \cdot
\Vert_{L^2 ,h_0}^2} = \log \, (q_0) - \log \, (q_1) &(14)
\cr 
= \log \det \left( \int_X e^{\vp} \, \lgl \a_i ,\a_j
\rgl \, \mu \right) + \log \det \left( \int_X e^{-\vp} \,
\lgl \b_i ,\b_j \rgl \, \mu \right) \, . \cr
}
$$
Proposition 1 follows from (5), (6) and (14).

\bigskip

\noindent {\bf 2.4.} {\bf Corollary 1.} {\it Under the assumption of
Proposition 1, if $b_1 =0$ and if the conjecture holds for one
choice of metric $h_X$ on $X$, it holds for any other choice of
metric on $X$.}

\bigskip

\noindent {\it Proof.} Assume $h_X$ gets replaced by $e^{\rho}
\, h_X$, where $\rho \in C^{\ify} (X)$. Then $\mu$ is replaced
by $\mu' = e^{\rho} \, \mu$. Let $(\a'_i)$ be an orthonormal
basis of $H^0 (X,L)$ for the $L^2$ metric defined by $h'_X$ and
$h_0$. We may write
$$
\a'_i = \sum_{j=1}^{b_0} a_{ij} \, \a_j \, ,
$$
$1\leq i \leq b_0$, where the square matrix $M = (a_{ij})$ is
independent of $\vp$. Therefore $A(\vp)$ is replaced by
$$
\eqalign{
B(\vp) = & \ {1\over 2} \int_X \vp \, dd^c (\vp) - \int_X \vp
(c_1 (T_X ,h_X) - dd^c \, \rho + c_1 (L,h_0)) \cr
+ & \ \log \det \left( \int_X e^{\vp +\rho} \, \lgl \a_i , \a_j
\rgl \, \mu \right)_{1 \leq i,j \leq b_0} \cr
+ & \ 2 \log \vert \det (M) \vert \, . \cr
}
$$
If we let $\psi = \vp + \rho$, since
$$
2 \int_X \vp \, dd^c (\rho) = \int_X (\vp \, dd^c \, \rho + \rho
\, dd^c \, \vp)
$$
we get
$$
B(\vp) = A(\psi) - {1\over 2} \int_X \rho \, dd^c (\rho) + 2
\log \vert \det (M)\vert \, .
$$
When $\rho$ is chosen, if $A(\psi)$ is bounded, so is $B(\vp)$.

\hfill q.e.d.

\vglue 1cm

\noindent {\bf 3. Proof of Theorem 1 in case i).}

\medskip

First notice that, because of (1),
 we can impose the condition
$$
\int_X \vp \, \mu = 0 \, . \leqno (15)
$$
On the other hand, we can choose the reference metric
$h_0$ in such a way that the form $c_1 (T_X ,h_X) + c_1
(C,h_0)$ is proportional to $\mu$. Together with (15),
this implies that the summand
$$
\int_X \vp \, (c_1 (T_X ,h_X)+c_1 (L,h_0))
$$
in Proposition 1 vanishes.

\medskip

Now let $A$ be an upper bound on $X$ for the $C^{\ify}$
functions $\vert \lgl \a_i ,\a_j \rgl \vert$, $1\leq i,j
\leq b_0$, and $\vert \lgl \b_i ,\b_j \rgl \vert$, $1\leq
i,j \leq b_1$. We get
$$
\log \det \left( \int_X e^{\vp} \, \lgl \a_i ,\a_j \rgl
\, \mu \right) \leq b_0 \, \log \int_X e^{\vp} \, \mu +
b_0 \, \log \, (A) + \log \, (b_0 !)
$$
and
$$
\log \det \left( \int_X e^{-\vp} \, \lgl \b_i ,\b_j \rgl
\, \mu \right) \leq b_1 \, \log \int_X e^{-\vp} \, \mu +
b_1 \, \log \, (A) + \log \, (b_1 !) \, .
$$
So Proposition 1 implies
$$
 \ A(\vp) 
\leq  \ {1\over 2} \int_X \vp \, dd^c \, \vp + b_0 \,
\log \int_X e^{\vp} \, \mu + b_1 \, \log \int_X e^{-\vp}
\, \mu + c_1 \, ,  \leqno (16) 
$$
for some constant $c_1 \geq 0$ independent of $\vp$.

\medskip

Let $\nb$ be the gradient defined by $h_X$. A local
computation shows that
$$
\int_X \vp \, dd^c (\vp) = -{1\over 4\pi} \int_X \vert \nb
\vp \vert^2 \, \mu \, . \leqno (17)
$$
We use now an inequality due to Fontana [3], Theorem 1.7,
which extends to arbitrary compact manifolds a result of
Moser and Trudinger for 
the sphere and open domains in $\Rb^n$ [7].
Namely, given any smooth real function $f$ on $X$ such
that
$$
\int_X \vert \nb f \vert^2 \, \mu \leq 1 \quad \hbox{and}
\quad \int_X f \, \mu =0 \, ,
$$
we have
$$
\log \int_X \exp (4\pi \, f^2) \, \mu \leq c_2 \, ,
\leqno (18)
$$
where $c_2$ is a constant which does not depend on $f$. From this
inequality it follows that, for any smooth real function
$g$ on $X$ such that $\int_X g \, \mu =0$,
$$
\log \int_X \exp (g) \, \mu \leq c_2 + {1\over 16\pi}
\int_X \vert \nb g \vert^2 \, \mu \, . \leqno (19)
$$
Indeed, if we let $B=\int_X \vert \nb g \vert^2 \, \mu$
and $f = g \, B^{-1/2}$, we have
$$
4\pi \, f^2 - g + {BÊ\over 16\pi} = 4\pi \left( f -
{\sqrt B \over 8\pi}\right)^2 \geq 0 \, .
$$
Therefore (18) gives
$$
\log \int_X \exp \, (g) \, \mu \leq \log \int_X \exp \,
(4\pi 
\, f^2) \, \mu + {B \over 16\pi} \leq c_2 + {1\over 16\pi}
\int_X \vert \nb g\vert^2 \, \mu \, .
$$
If we apply the inequality (19) to $\vp$ and $-\vp$ we
get, from (16) and (17), the inequality 
$$
\eqalign{
 \ A(\vp) 
\leq  \ -{1\over 8\pi} \int_X \vert \nb \vp \vert^2 \, 
\mu + {b_0 \over 16\pi} \int_X \vert \nb \vp \vert^2 \,
\mu + {b_1 \over 16\pi} \int_X \vert \nb \vp \vert^2 \,
\mu 
+  \ c_1 + 2 \, c_2 \, . \cr
}
$$
When $b_0 + b_1 \leq 2$ we conclude that
$$
A(\vp) \leq c_1 +
2 \, c_2 \, .
$$

\hfill q.e.d.

\vglue 1cm

\noindent {\bf 4. Proof of Theorem 1 in case ii).}

\medskip

\noindent {\bf 4.1.} We assume that $X=\Pb^1$ is the
complex projective line and that $L=\Oc (n)$, $n\geq 1$.
Then $b_1 =0$ and $b_0 = n+1$. Furthermore
$$
H^0 (X,\Oc (n)) = S^n \, H^0 (X,\Oc (1))
$$
is the space of homogeneous polynomials of degree $n$ in
two variables. Consider the canonical exact sequence of
sheaves
$$
0 \ra \Oc (-1) \ra \Cb^2 \ra \Oc (1) \ra 0 \, ,
$$
and denote by $A,B \in H^0 (X,\Oc (1))$ the images of the
vectors $(1,0)$, $(0,1)$ in $\Cb^2 = H^0 (X,\Cb^2)$.
Choose on $\Oc (1)$ the metric $h_0$ induced by the
standard metric on $\Cb^2$. At a point $P$ with
homogeneous coordinates $(u,v)$ in $\Pb^1$ the lift of
$A$ to $\Cb^2$ which is orthogonal to the vector $(u,v)
\in
\Oc (-1)_P$ is given by
$$
A^{\perp} = (\vert v \vert^2 , -\ov u \, v) \ (\vert u
\vert^2 + \vert v \vert^2)^{-1} \, .
$$
Similarly $B$ lifts to
$$
B^{\perp} = (-u \, \ov v , \vert u \vert^2) \
(\vert u \vert^2 + \vert v \vert^2)^{-1} \, .
$$
If $P$ lies on the affine line $\Ab^1 \sbs \Pb^1$ with
affine coordinate $z$, these vectors become
$$
A^{\perp} = (1,-\ov z ) \, N^{-1}
$$
and
$$
B^{\perp} = (-z ,\vert z \vert^2) \, N^{-1} \, ,
$$
where $N=\vert z \vert^2 +1$. The scalar products of the
sections $A$ and $B$ of $H^0 (X,\Oc (1))$ at the point
$P$ are thus given by
$$
\ldisplaylinesno{
\lgl A,A \rgl = A^{\perp} \cdot A^{\perp} = N^{-1} \, ,
&(20) \cr
\lgl A,B \rgl = A^{\perp} \cdot B^{\perp} = -\ov z
\, N^{-1} \, , \cr
}
$$
and
$$
\lgl B,B \rgl = B^{\perp} \cdot B^{\perp} = \vert z
\vert^2 \,  N^{-1} \, .
$$

\medskip

An orthogonal basis of $H^0 (X,\Oc (n))$ is the set of
monomials $(A^i \, B^j , i+j =n)$, where $A^i \, B^j$ is
the symmetrization of the vector $A^{\ot i} \ot B^{\ot
j}$ in $H^0 (X,\Oc (1))^{\ot n}$. Using (20) we see that
the standard metric on $S^n \, H^0 (X,\Oc (1))$, which is
(a constant multiple of) the $L^2$-metric on $H^0 (X,\Oc
(n))$, is such that
$$
\lgl A^i \, B^j , A^k \, B^{\ell} \rgl = \left({n \atop
j}\right) \, \left({n \atop \ell}\right) \, {(-z)^j \,
(-\ov z )^{\ell} \over N^n} \, . \leqno (21)
$$

\bigskip

\noindent {\bf 4.2.} To prove Theorem 1 ii), we may assume $n >
0$, because of (2). By the argument of Corollary 1, we may also
assume that both $h_0$ and $h_X$ are the standard metrics. In
particular
$$
\mu = {d z \, d \bar z \over 2i\pi (1+\vert z \vert^2)^2}
$$
and $c_1 (L,h_0)$ is a multiple of $\mu$.
\noindent By 1.2. a) we may finally assume that
$$
\int_X \vp \, \mu =0 \, .
$$
From the previous section, an orthonormal basis $(\a_i)$ of $H^0
(X,L)$ is given by the elements
$$
A^k \, B^{\ell} / \Vert A^k \, B^{\ell} \Vert_{L^2} \ , \quad
k+\ell =n \, .
$$
From (21), since $\vp$ and $\mu$ are invariant under rotation,
we conclude that
$$
\int_X e^{\vp} \, \lgl \a_i ,\a_j \rgl \, \mu = 0
$$
when $i\not= j$. From Proposition 1 we get
$$
\leqalignno{
A(\vp) = & \ {1\over 2} \int_X \vp \, dd^c (\vp) &(22) \cr
+ & \ \sum_{i=0}^{n} \log \int_X e^{\vp} \, {\vert z \vert^{2i}
\over (1+\vert z \vert^2)^{n+2}} \, dz \, d\bar z 
+  \ c_1 \cr
}
$$
where $c_1 ,c_2$ etc $\ldots$ will denote constants independent
of $\vp$. If we take polar coordinates $z = r \, e^{i\t}$ and if
we make the change of coordinates $r = e^{t/2}$, $t\in \Rb$, we
may write
$$
\matrix{
\hfill \vp (x) &= & f(t) \, , \hfill \cr
\cr
\hfill \rho (t) &= & (e^{t/2} + e^{-t/2})^{-2} \, , \hfill \cr
\cr
\hfill \rho_i (t) &= &{e^{it} \over (1+e^t)^n} \,
\rho (t) \, , \hfill \cr
} \leqno (23)
$$
in which case
$$
\int_X \vp \, dd^c (\vp) = - \int_{-\ify}^{+\ify} \dot f (t)^2 \,
dt
$$
(where $\dot f (t)$ is the derivative of $f(t)$) and
$$
\int_X e^{\vp} \, {\vert z \vert^{2i} \over (1+\vert
z\vert^2)^{n+2}} \, dz \, d\bar z = 2\pi \int_{-\ify}^{+\ify}
e^{f(t)} \, \rho_i \, (t) \, dt \, .
$$
We conclude that
$$
\leqalignno{
A(\vp) = & \ -{1\over 2} \int_{-\ify}^{+\ify} \dot f (t)^2 \, dt
&(24) \cr
+ & \ \sum_{i=0}^n \log \int_{-\ify}^{+\ify} e^{f(t)} \, \rho_i
\, (t) \, dt 
+  \ c_2 \, . \cr
}
$$
Furthermore
$$
\int_X \vp \, \mu = \int_{-\ify}^{+\ify} f(t) \, \rho \, (t) \,
dt =0 \, . \leqno (25)
$$

\medskip

\noindent {\bf 4.3.} Let $A>0$ be such that
$$
A^2 = \int_{-\ify}^{\ify} \dot f (t)^2 \, dt \, . \leqno (26)
$$
We first deduce from (25) that there is a constant $c_3$ such
that
$$
|f(0)| \leq c_3  A \leqno (27)
$$
(compare [7] (8)). Indeed, the Cauchy-Schwarz inequality implies
$$
\left\vert \int_s^t \dot f (t) \, dt \right\vert^2 \leq \int_s^t
\dot f (t)^2 \, dt \int_s^t dt
$$
i.e.
$$
(f(t)-f(s))^2 \leq A^2 \, \vert t-s \vert
$$
for all $s$ and $t$, hence
$$
-A \, \sqrt{\vert t-s \vert} \leq f(t) - f(s) \leq A \,
\sqrt{\vert t-s \vert} \, . \leqno (28)
$$
We multiply these inequalities by $\rho (s)$ and we integrate
with respect to $s$. Since $\int_{-\ify}^{\ify} \rho (s) \, ds
=1$ we get from (25) that
$$
\vert f(t) \vert \leq A \int_{-\ify}^{+\ify} 
\sqrt{\vert t-s \vert} \,
\rho \, (s) \, ds
$$
and (27) follows when $t=0$.

\bigskip

\noindent {\bf 4.4.} {\bf Lemma 2.}{\it There exists a function
$u(t)$, $t\in \Rb$, $t\geq 0$, which is $C^1$ and such that

\smallskip

\item{i)} $u(0) = f(0)$, $u(+\ify ) = f(+\ify )$;

\smallskip

\item{ii)} $\dot u (t) \geq 0$, $\dot u (t)$ is nonincreasing;

\smallskip

\item{iii)} $u(t) \geq f(t)$;

\smallskip

\item{iv)} $\int_0^{\ify} \dot u (t)^2 \, dt = \int_0^{\ify}
\dot f (t)^2 \, dt$.}

\bigskip

\noindent {\it Proof of Lemma 2.} Let ${\dot f}^*$ be the
nonincreasing rearrangement of $\dot f$ on $[0,+\ify [$ (cf.
e.g. [8]). In other words ${\dot f}^*$ is the nonincreasing
function on $[0,+\ify [$ such that, for all $y \geq 0$, $({\dot
f}^*)^{-1} \, (y)$ is the Lebesgue measure of the set of numbers
$t\in [0,+\ify [$ such that $\dot f (t) \geq y$. Since $\dot f$
is continuous, the same is true for ${\dot f}^*$ and we may
define
$$
u(t) = f(0) + \int_0^t {\dot f}^* \, (s) \, ds \, . \leqno (29)
$$
It is clear that $u(0) = f(0)$, and 
the standard equalities ([8],
Lemma 2.2)
$$
\int_0^{+\ify} ({\dot f}^*)^k \, (s) \, ds = \int_0^{+\ify}
{\dot f}^k \, (s) \, ds \, ,
$$
for $k=1,2$, imply that $u(\ify) = f(\ify)$ and that iv) holds.
Property ii) is a consequence of the definitions and iii) is
equivalent to
$$
\int_0^t {\dot f}^* \, (s) \, ds \geq \int_0^t \dot f (s) \, ds
\, ,
$$
a well-known property of rearrangements. 

\hfill q.e.d.

\bigskip

\noindent {\bf 4.5.} To bound the quantity $A(\vp)$ in (24) we
may now assume, by Lemma 2, that $\dot f (t) \geq 0$, that $\dot
f (t)$ is a nondecreasing function when $t \leq 0$ and a
nonincreasing function when $t \geq 0$. Note that
$$
\log (a+b) \leq \log (2) + \log^{+} (a) + \log^{+} (b) \, ,
$$
where $\log^{+} = Max(log,0)$.
Therefore, Lemma 3 below, when applied to $f(t)$ and $f(-t)$,
gives
$$
\eqalign{
{} & \ \sum_{i=0}^{n} \log \int_{-\ify}^{+\ify} e^{f(t)} \,
\rho_i \, (t) \, dt \cr
\leq & \ \sum_{i=0}^{n} 
\log^{+} \int_0^{+\ify} e^{f(t)} \, \rho_i
\, (t) \, dt \cr
+ & \ 
\sum_{i=0}^{n} \log^{+} \int_{-\ify}^0 e^{f(t)} \,
\rho_i \, (t) \, dt + (n+1) \, \log (2) \cr
\leq & \ 
\sum_{i=0}^{n} \log^{+} \int_0^{+\ify} e^{f(t)} \,
e^{-(i+1)t} \, dt \cr
+ & \ 
\sum_{i=0}^{n} \log^{+} \int_{-\ify}^0 e^{f(t)}
\, e^{(n-i-1)t} \, dt + (n+1) \, \log(2) \cr
\leq & \ 2(n+1) \, \vert f(0)\vert + \left( {1\over 2} - {1\over
70n^2} \right) \, A^2 + c_4 \, . \cr
}
$$
From (27) we conclude that
$$
2(n+1) \, \vert f(0)\vert \leq 2(n+1) \, c_3 \, A \leq {1\over
70n^2} \, A^2 + c_5 \, .
$$
Therefore we get
$$
\sum_{i=0}^{n} \log \int_{-\ify}^{+\ify} e^{f(t)} \, \rho_i \,
(t) \leq {1\over 2} \, A^2 + c_6 \, ,
$$
i.e. (by (24) and (26)) $A(\vp)$ is bounded from above and
Theorem 1 ii) is proved.

\bigskip

\noindent {\bf 4.6.} {\bf Lemma 3.} {\it Let $M \geq 1$ be an integer
and let $u : \Rb^+ \ra \Rb$ be a $C^1$ map such that $\dot u$ is
$L^2$ and nonincreasing. Define
$$
X = \sum_{j=0}^{M} \log \int_0^{+\ify} \exp (u(t) - (j+1)t)
\, dt \, .
$$
Then there exists a constant $C \geq 0$ such that
$$
X \leq (M+1) \, \vert u(0)\vert + \left( {1\over 2} - {1\over
70M^2} \right) \int_0^{+\ify} \dot u (t)^2 \, dt + C \, . \leqno
(30)
$$}

\bigskip

\noindent {\it Proof of Lemma 3.} For any integer $k \geq 1$ we
let
$$
\lb_k = 1+{1\over 5k^2}
$$
and
$$
\mu_k = 1-{1\over 4k} \, .
$$
Note that $\lb_k\cdot  k + \mu_k > k$ and
 $r_k := k+1 - \lb_k\cdot k 
-\mu_k = {1\over 4k} - {1\over 5k} > 0$.

\smallskip

\noindent Define $N \geq 0$ as the smallest integer such that
$$
\dot u (0) \leq \lb_{N+1} \cdot (N+1) + \mu_{N+1} \, .
$$

\medskip

\noindent ${\ast}$ 
 If $j \geq N+1$ and $t \geq 0$ we have
$$
\dot u (t) - (j+1) \leq \dot u (0) - (j+1) \leq \dot u (0) -
(N+2) \leq -r_{N+1} \, .
$$
Therefore
$$
\log \int_0^{\ify} \exp (u(t)-(j+1) \, t) \, dt \leq u(0) + \log
\int_0^{\ify} \exp (-r_{N+1} \, t) \, dt = u(0) + c \, .
$$

\medskip

\noindent  ${\ast}$ If $N =0$ and $0 \leq j < N+1$, i.e. $j=0$,
 we know from (28) that
$$u(t) \leq u(0) + \sqrt{I\cdot t}$$
where
$$
I = \int_0^{\ify} \dot u (t)^2 \, dt \, .
$$
Therefore, by completing a square we get
$$u(t)-t \leq u(0) + \sqrt{I\cdot t} -t \leq u(0) + {3I \over 8}
- {t \over 3}, $$
from which it follows that
$$
\log \int_0^{\ify} \exp (u(t)-t) \, dt \leq 
u(0) + {3I \over 8} +
\log \int_0^{\ify} \exp (-t/3) \, dt
$$
hence, if $N=0$,
$$X \leq (M+1)u(0)+{3I\over8}+C\  .$$ 

\medskip

\noindent  ${\ast}$ If $N \geq 1$, we let $\lb = \lb_N$, $\mu = \mu_N$
and we choose real numbers $x_0 \geq x_1 \geq \cdots \geq x_N >
0$ such that
$$
\dot u (x_j) = \lb\cdot j + \mu \ , \quad 0 \leq j < N \, . \leqno (31)
$$
Then if $0 \leq j < N$, we have
$$u(t)-(j+1)t = u(t)-(\lb\cdot j + \mu)t + 
(\lb\cdot j + \mu - (j+1))t$$
$$
\leq  u(x_j)-(\lb\cdot j + \mu)x_j + 
(\lb\cdot N + \mu - (N+1))t
$$
hence
$$
\eqalign{
\log \int_0^{\ify} \exp (u(t)-(j+1) \, t) \, dt \leq & \cr
u(x_j)-(\lb\cdot j + \mu) \, x_j 
+ & \ \log \int_0^{\ify} \exp (-r_N \, t) \, dt \cr
= & \vp (x_j) + c'\, ,
\cr 
}
$$
where 
$$
\vp (x) = u(x) - \dot u (x) \, x \, .
$$
Therefore
$$
X \leq (M-N) \, (u(0)+c) + Nc' +Y
$$
where
$$
Y = \sum_{j=0}^N \vp (x_j) \, .
$$

Using (31) we can write
$$
\eqalign{
Y = & \ \sum_{j=0}^N \vp (x_j) = {1\over \lb} \sum_{j=1}^N \vp
(x_j) \, (\dot u (x_j) - \dot u (x_{j-1})) + \vp (x_0) \cr
= & \ {1\over \lb} \sum_{j=1}^{N-1} \dot u (x_j) \, (\vp (x_j) -
\vp  (x_{j+1})) - {1\over \lb} \, \vp (x_1) \, \dot u (x_0) \cr
+ & \ {1\over \lb} \, \dot u (x_N) \, \vp (x_N) + \vp (x_0) \cr
= & \ {1\over \lb} \sum_{j=0}^{N-1} \dot u (x_j) \, (\vp (x_j) -
\vp (x_{j+1})) + {1\over \lb} \, \dot u (x_N) \, (\vp (x_N) -
\vp(0)) \cr
+ & \ \left( N + {\mu \over \lb}\right) \, u(0) + \vp (x_0) \,
\left( 1-{\dot u (x_0) \over \lb} \right) \cr
\leq & \ {1\over \lb} \int_0^{x_0} \dot u \, \dot{\vp} \, dt +
\left( N + {\mu \over \lb}\right) \, u(0) + \vp (x_0) \, \left(
1-{\dot u (x_0) \over \lb} \right) \, . \cr
}
$$
But
$$
\eqalign{
\int_0^{x_0} \dot u \, \dot{\vp} \, dt = & \ -\int_0^{x_0}
t \dot u
 \ddot u  \, dt \cr
= & \ \left[ -{1\over 2} t \dot u \dot u  \right]_0^{x_0} +
{1\over 2} \int_0^{x_0} \dot u (t)^2 \, dt \cr
= & \ -{1\over 2} \, \dot u (x_0)^2 \, x_0 + {1\over 2}
\int_0^{x_0} \dot u (t)^2 \, dt \, . \cr
}
$$
So we get
$$
Y \leq {1\over 2\lb} \, I + \left( N + {\mu \over \lb} \right)
\, u(0) + R \leqno (32)
$$
where
$$
\eqalign{
R = & \ -{1\over 2\lb} \, \dot u (x_0)^2 \, x_0 + \vp (x_0) \,
\left( 1- {\dot u (x_0) \over \lb} \right) \cr
= & \ \left( 1 -{\mu \over \lb} \right) \, u (x_0) + \left(
{\mu^2 \over 2\lb} -\mu \right) \, x_0 \, . \cr
}
$$
Now, by (28),
$$
u(x_0) \leq u(0) + \sqrt{x_0\cdot I}\, ,
$$
and, by completing the square,
$$
\a \, \sqrt{x_0} - \b \, x_0 \leq {\a^2 \over 4\b} \qquad
\hbox{for all} \ \a ,\b > 0 \, .
$$
Therefore
$$
R \leq \left( 1-{\mu \over \lb} \right) \, u(0) + {(1-\mu /
\lb)^2 \over 4 \left( \mu - {\mu^2 \over 2\lb}\right)} \, I \, .
$$
Using (32) we get
$$
Y \leq (N+1) \, u(0) + A \cdot I
$$
where
$$
A = {1\over 2\lb} + {\left( 1-{\mu \over \lb} \right)^2 \over 4
\left( \mu - {\mu^2 \over 2\lb} \right)} \, .
$$
From the values of $\lb = \lb_N$ and $\mu = \mu_N$ we compute
$$
A \leq {1 \over 2} - {1\over 70N^2} \leq {1\over 2} - {1\over
70M^2} \, .
$$
Therefore
$$
X \leq (M+1) \, u(0) + \left( {1\over 2} - {1\over 70M^2}
\right) \int_0^{\ify} \dot u (t)^2 \, dt + C \, .
$$
\hfill q.e.d.

\vglue 1cm

\noindent {\bf 5. Flat bundles}

\medskip

According to Bismut and Zhang [2]  a flat
$C^{\ify}$ bundle $(E, \nb)$,  $\nb^2 =0$, 
 together with a $C^{\ify}$ metric $h$ on $E$
on a $C^{\ify}$ manifold $M$ 
is the analog in the differentiable
category of a hermitian holomorphic bundle on  a complex
manifold.
 Inequalities similar to our conjecture might also hold
in this case, but in some cases they must be lower bounds rather 
than upper bounds, as the following example suggests.

\smallskip

Let $M=S^1$ be the circle, and let $E=\Cb$ be the trivial line
bundle on $M$. We equip $M$ with its standard metric and $E$
with an arbitrary metric $h$. The connection $\nb = d$ has an
adjoint $d^*$ (depending on $h$), and we consider the Laplace
operator $\D = d^* \, d$ on $C^{\ify} (M)$, and its regularized
determinant $\det' (\D)$.

\bigskip

\noindent {\bf Proposition 3.}{\it There is a constant $C(E)$ such
that, for any choice of a metric $h$ on $E = \Cb$,
 the following inequality holds	
$$
{\textstyle \det'} (\D) \leq C(E) \, .
$$}

\bigskip

\noindent {\it Proof.} We have $H^0 (S^1 ,\Cb) \sm H^1 (S^1
,\Cb) \sm \Cb$. On $\lb (E) = H^0 (S^1 ,\Cb) \ot H^1 (S^1
,\Cb)^*$ we define the Quillen metric by
$$
\Vert \cdot \Vert_Q^2 = \Vert \cdot \Vert_{L_2}^2 \,
{\textstyle \det'} (\D)^{-1}
$$
as in \S 1 above [2].

\smallskip

\noindent If $\vp \in C^{\ify} (M)$ we let $h$ be the metric on
$E = \Cb$ such that $h(1,1) = \exp (\vp)$, and we denote by
$\D_{\vp}$ the corresponding Laplace operator on $C^{\ify}
(M,\Cb)$. We define
$$
\eqalign{
A(\vp) = & \ \log ({\textstyle \det'} (\D_{\vp})) - \log
({\textstyle \det'} (\D_0)) \cr 
= & \ \log {\Vert \cdot \Vert_{L^2,\vp}^2 \over \Vert
\cdot \Vert_{L^2 ,0}^2} - \log \, {\Vert \cdot \Vert_{Q,\vp}^2
\over \Vert \cdot \Vert_{Q,0}^2} \, . \cr
}
$$
According to [2] Theorem 0.1 we have
$$
\log \, {\Vert \cdot \Vert_{Q,\vp} \over \Vert \cdot
\Vert_{Q,0}} = 0
$$
and a computation similar to \S 2.3 gives
$$
\log \, {\Vert \cdot \Vert_{L^2 ,\vp}^2 \over \Vert \cdot
\Vert_{L^2 ,0}^2} = \log \int_{S^1} e^{\vp (x)} \, dx + \log
\int_{S^1} e^{-\vp (x)} \, dx \, ,
$$
where $dx$ is the Haar measure of length one.
The Cauchy-Schwarz inequality implies
$$
\left( \int_{S^1} dx \right)^2 \leq \left( \int_{S_1} e^{\vp
(x)} \, dx \right) \, \left( \int_{S^1} e^{-\vp (x)} \, dx
\right)
$$
from which we conclude that $A(\vp) \geq 0$.

\hfill q.e.d.

\bigskip

\noindent {\bf References.}

\medskip

\item{[1]} Bismut J.-M.,  Gillet H., Soul\'e C.:
Analytic torsion and holomorphic determinant bundles
I. Bott-Chern forms and analytic torsion,
 {\it Comm. in Math. Physics
} {\bf 115} (1988), 49-78.

\item{[2]} Bismut J.-M., Zhang W.: An extension
of the Cheeger-Muller theorem,
{\it Ast\'erisque} {\bf 205} (1992),
Soc. Math. de France, Paris.

\item{[3]} Fontana L.: Sharp borderline Sobolev
inequalities on compact Riemannian manifolds, {\it Comm.
Math. Helv.} {\bf 68} (1993), 415-454.

\item{[4]} Gillet H., Soul\'e C.:
Characteristic classes for algebraic vector
bundles with hermitian metrics I, 
{\it Annals of Math.} {\bf 131} (1990), 163-203.

\item{[5]} Gillet H., Soul\'e C.: An arithmetic Riemann-Roch
 theorem, 
 {\it Inventiones Math.} {\bf 110} (1992), 473-543.

\item{[6]} Kiessling M.K.H.: 
Statistical mechanics of classical particles with
logarithmic interactions,
 {\it Comm. Pure Appl. Math.} {\bf 46} (1993), 27-56.

\item{[7]} Moser J.: A sharp form of an inequality by N.
Trudinger, {\it Indiana Math. J.} {\bf 20} (1971),
1077-1092.

\item{[8]} O'Neil R.: Convolution operators and $L(p,q)$
spaces, {\it Duke Math. J.} {\bf 30} (1963), 129-142.

\item{[9]} Quillen D.: Determinants of Cauchy-Riemann
operators over a Riemann surface, {\it Funct. Anal.
Appl.} (1985), 31-34.

\item{[10]} Soul\'e C.: A vanishing theorem on arithmetic surfaces,
 {\it Inventiones Math.} {\bf 116} (1994), 577-599.

\bigskip

H.G.:  Department of Mathematics, Statistics, and Computer Science,
 University of Illinois at Chicago, 851 S. Morgan Street,
Chicago, IL 60607-7045, U.S.A.

\medskip

C.S. : CNRS, Institut des Hautes \'{E}tudes Scientifiques,  35, Route de
Chartres, 91440, Bures-sur-Yvette, France

\bye